\documentclass{article}

\usepackage{amsmath,amsfonts,amsthm}

\numberwithin{equation}{section}

\theoremstyle{plain}
\newtheorem{lemma}[equation]{Lemma}

\theoremstyle{definition}
\newtheorem{definition}[equation]{Definition}

\title{Points on Hemispheres}
\author{Jan Fricke}

\begin{document}

\maketitle

\begin{abstract}
  We will show that for any $n\ge N$ points on the $N$-dimensional sphere
  $S^N$ there is a closed hemisphere which contains at least
  $\lfloor\frac{n+N+1}{2}\rfloor$ of these points. This bound is sharp
  and we will calculate the amount of sets which realize this value.

  If we change to open hemispheres things will be easier.
  For any $n$ points on the sphere there is an open hemisphere which
  contains at least $\lfloor\frac{n+1}{2}\rfloor$ of these points,
  independent of the dimension. This bound is sharp.
\end{abstract}

\section{Introduction}
In problem 451342 of the German Mathematical Olympiad one has to
calculate the maximum of the minimal pairwise distance of five points
on a (two-dimensional) sphere. For proving the result it was helpful
to use the property that there is a (closed) hemisphere containing
at least four of these points.

We now want to generalize this property to all numbers of points and
all dimensions. The \emph{$N$-dimensional sphere} $S^N$ is the set of
points in $\mathbb{R}^{N+1}$ with (Euclidean) distance $1$ from the
origin. Any (hyper-)plane through the origin splits the sphere into
two \emph{hemispheres}. The closed hemispheres contain the
intersection set, which is called \emph{great circle}, while the open
hemispheres do not contain the great circle.

Any hemisphere has a unique \emph{pole} $p$. Then the points $x$ of
the great circle, the closed and open hemisphere can be characterized
by $\langle x,p \rangle = 0$, $\langle x,p \rangle \ge 0$,
and $\langle x,p \rangle > 0$, respectively.

\section{The lower bound}
\begin{lemma}
  For any $n<N$ points on the $N$-dimensional sphere $S^N$, 
  there is a closed hemisphere which contains all of these points.
\end{lemma}

\begin{proof}
  For $n<N$ points on the $N$-dimensional sphere there is always a
  great circle containing all these points.
\end{proof}

\begin{lemma}
  For any $n$ points on the $N$-dimensional sphere
  $S^N$, $n\ge N$, there is a closed hemisphere which contains at least
  $\lfloor \frac{n+N+1}{2} \rfloor$ of these points.
\end{lemma}

\begin{proof}
  Choose any $N$ of these points. Then there is a plane through these
  points and the origin. So there is a great circle containing
  these $N$ points. The remaining $n-N$ points distribute to the two
  hemispheres, hence the is a hemisphere containing at least
  $\lfloor \frac{n-N+1}{2} \rfloor$ of them. Together with the initial
  $N$ points there are at least $\lfloor \frac{n+N+1}{2} \rfloor$ on
  this hemisphere.
\end{proof}

\section{The sharpness of the bound}
Now we will check that the bound from the previous section is
sharp.

Take any great circle and move it continuously until it contains $N$
points. Then the number of points in any of the two hemispheres is
non-decreasing. So no hemisphere contains more than
$\lfloor \frac{n+N+1}{2} \rfloor$ points if and only if all the
hemispheres whose border contains $N$ points have this property.

\begin{definition}
  A set of $n$ points on the $N$-dimensional sphere $S^N$ is called
  \emph{equator-balanced set}, if any great circle through $N$ of
  these points cuts the remaining points into halves. (If the number
  is odd then the two sets differ by one.)
\end{definition}

Hence we have to show the existence of such sets.

\begin{lemma}\label{lem:exist}
  For any $n>N$ there is an equator-balanced set of $n$ points on the
  $N$-sphere.
\end{lemma}

\begin{proof}
  One can specify such a set. Let $y_i=(-1)^i\cdot(1,i,i^2,\dotsc,i^N)$,
  $i=1,\dotsc,n$ be $n$ vectors in $\mathbb{R}^{N+1}$, and $x_i$ the
  corresponding normal vectors on the $S^N$.

  Fix any $N$ of these points. In which hemisphere the point $x_i$
  lies is determined by the determinant formed by $x_i$ and the $N$
  points. Apart from some factors this is a Vandermonde determinant.

  Replacing $x_i$ by the next remaining point $x_{i+j}$ the signum of
  the factor changes by $(-1)^j$ and in the Vandermonde determinant
  $j-1$ factors change their sign. So the determinant has the
  opposite sign.

  Hence the remaining points lie alternating in the two hemispheres,
  and the set is equator-balanced.
\end{proof}

\section{The density of equator-balanced sets}
Now we can ask with which probability a set is equator-balanced.

\begin{definition}
  By $p(N,n)$ we denote the probability that a set of $n$ independent
  uniformly distributed points on the $N$-sphere is equator-balanced.
\end{definition}

\begin{lemma}
  $p(n,N)$ exists and is always positive.
\end{lemma}

\begin{proof}
  In the configuration space of $n$-sets on the $N$-sphere let
  $M(n,N,k)$ be the set of all sets where any closed hemisphere
  contains at most $k$ points. This is equivalent to: any open
  hemisphere contains at least $n-k$ points. But this is an open
  condition, i.e. if a configuration is slightly changed then it also
  fulfills this condition. Hence $M(n,N,k)$ is an open set.

  Since the equator-balanced sets constitute the set $M(n,N,k)$ with
  $k=\lfloor \frac{n+N+1}{2} \rfloor$ and by lemma \ref{lem:exist}
  this set is non-empty, it is measurable with positive
  measurable. Hence $p(n,N)$ exists and is always positive.
\end{proof}

In some special cases $p(N,n)$ can be calculated. Obviously,
if $n\le N+1$ then all sets are equator-balanced. So we will restrict
to the case $n\ge N+2$.

\begin{lemma}
  $p(N,N+2)=2^{-(N+1)}$.
\end{lemma}

\begin{proof}
  A set of $N+2$ points on the $N$-sphere is not equator-balanced if
  and only if all the points lie in a common hemisphere.

  Fix $N+1$ of the points. They constitute a spherical
  $(N+1)$-simplex. The only possibility for the last point not lying
  in a common hemisphere with the other points is to lie in the
  $(N+1)$-simplex antipodal to the simplex.

  Hence the probability $p(N,N+2)$ is the expectation value of the
  volume of a random spherical $(N+1)$-simplex relative to the volume
  of the $N$-sphere. For independent random points
  $x_1,\dotsc,x_{N+1}$ we have
  \begin{equation}
    E(V(x_1,\dotsc,x_{N+1})) =  E(V(\pm x_1,\dotsc,\pm x_{N+1}))
  \end{equation}
  for any combination of plus and minus signs. But all these $2^{N+1}$
  simplices assembled together always gives the whole $N$-sphere, so
  the sum of these (equal) $2^{N+1}$ expectation values is $1$, thus
  $p(N,N+2)=E(V(x_1,\dotsc,x_{N+1}))=2^{-(N+1)}$.
\end{proof}

\begin{lemma}
  $p(1,1+2k)=4^{-k}$.
\end{lemma}

\begin{proof}
  Any of the $1+2k$ can be replaced by its antipodal point. Consider
  all $2^{1+2k}$ combinations where each point can be replaced by its
  antipodal point. For all these combinations the probability to be
  equator-balanced is the same. We will show, that independently from
  the original configuration there are exactly $2$ equator-balanced
  configurations among them. So the probability is $2/2^{1+2k} = 4^{-k}$.

  Fix any diameter of the circle and count the points in one
  semicircle. Now rotate the diameter continuously half around and
  notice the change of the points in the semicircle. This numbers give
  a sequence of $2+2k$ numbers, where successive numbers differ by $1$
  and the first and last number add up to $1+2k$. The other way round
  any such sequence gives a configuration. The configuration is
  equator-balanced if and only if the sequence contains only the
  numbers $k$ and $k+1$. There are exactly two such sequences:
  $k,k+1,\dotsc,k+1$ and $k+1,k,\dots,k$.
\end{proof}

In the preceeding examples $n-N$ was even, i.e. we coped with
\emph{real} equator-balanced sets. For the next interesting examples
of this kind a stochastical simulation gave the results in table \ref{table:stoch}.

\begin{table}[!h]
  \centering
  \begin{tabular}{|c|c|r|r|c|c|}
    \hline
    $N$ & $n$ & trials & success & $1/p(N,n)$ & precision \\
    \hline
    2 & 6 & 287951134242 & 1708252518 & 168.5647 & 0.012 \\
    2 & 8 & 293892632084 &   23669718 & 12416.39 & 7.65 \\
    3 & 7 & 115638779856 &   42369783 & 2729.27  & 1.25 \\
    \hline
    2 & 5 & 889631743 & 277996246 & 3.20015 & 0.00057 \\
    2 & 7 & 3558495944 & 254538093 & 13.98021 & 0.0026 \\
    2 & 9 & 11535004949 & 127320713 & 90.598 & 0.024 \\
    \hline
  \end{tabular}
  \caption{Stochastic simulation for $p(N,n)$. Here the precision is
    the estimated $3\sigma$-value for $1/p(N,n)$.}
  \label{table:stoch}
\end{table}

If $n-N$ is odd, then the amount of points in the two hemispheres
bounded by a great circle trough $N$ points differ by one. Hence there
are much more possibilities for such configurations and the
probability is much higher. But only in the case $N=1$ we were able to
determine the exact probability.

\begin{lemma}
  $p(1,2+2k)=2^{-k}$.
\end{lemma}

\begin{proof}
  As in the previous proof the calculation can be reduced to counting
  sequences. We have now sequences of length $3+2k$, where first and
  last number add up to $2+2k$. A configuration is equator-balanced if
  and only if the sequence contains only the numbers $k$, $k+1$ and
  $k+2$.

  Hence any second number in the sequence is $k+1$, and the remaining
  numbers can be choosen arbitrarily $k$ or $k+2$ (except the last
  number). We have two possibilities to choose the ``$k+1$-numbers''
  and again $2^{k+1}$ possibilities to choose the remaining numbers.

  So we have $2^{k+2}$ such sequences, and the probability is
  $2^{k+2}/2^{2+2k}=2^{-k}$.
\end{proof}

\section{Open hemispheres}
For open hemispheres the problem is much easier, because it is
possible to ``hide'' points by pairs of antipodal points.

\begin{lemma}\label{lem:nointersection}
  For any finite set of points on a sphere there is a great circle
  that does not contain any of these points.
\end{lemma}

\begin{proof}
  The poles of a great circle containing a given point lie on a great
  circle having the given point as a pole. Hence the poles of all
  great circles containing at least one of the points lie on a union of
  a finite number of great circles. So there are poles left for that
  the corresponding great circle does not contain any of these points.
\end{proof}

\begin{lemma}
  For any $n$ points on the sphere, there is an open hemisphere which
  contains at least $\lfloor \frac{n+1}{2} \rfloor$ of these points.
  This bound is best possible.
\end{lemma}

\begin{proof}
  By Lemma \ref{lem:nointersection} there is a great circle containing
  none of these points. So one of the corresponding hemispheres
  contains at least $\lfloor \frac{n+1}{2} \rfloor$ of these points.

  If the points lie pairwise antipodal (and one additional point if
  $n$ is odd), then any open hemisphere contains at most one point of
  any pair. Hence in this case any open hemisphere contains at most
  $\lfloor \frac{n+1}{2} \rfloor$ of these points.
\end{proof}

\end{document}